\documentclass[12pt]{article}

\usepackage{amsmath}
\usepackage{amsthm}
\usepackage{amsfonts}

\newenvironment{myproof}{\noindent {\it Proof} }{$\Box$}
\newtheorem{theorem}{Theorem}[section]
\newtheorem{lemma}{Lemma}[section]

\theoremstyle{remark}
\newtheorem{rem}{Remark}[section]
\newtheorem{corr}{Corollary}[section]

\newcommand{\suk}{\underset{k\geq 0}{\sup}}
\newcommand{\suin}{\underset{n\geq 0}{\sup}}
\newcommand{\sun}{\sum_{n=1}^{\infty}}
\newcommand{\su}{\sum_{n=0}^{\infty}}
\newcommand{\zi}{Z_{\infty}}

\newcommand{\suik}{\underset{k\geq 1}{\sup}}
\newcommand{\lin}{\underset{n\rightarrow\infty}{\lim}}
\newcommand{\lix}{\underset{x\rightarrow\infty}{\lim}}
\newcommand{\lisx}{\underset{x\to\infty}{\lim\sup}}
\newcommand{\liix}{\underset{x\to\infty}{\lim\inf}}
\newcommand{\mmp}{\mathbb{P}}
\newcommand{\me}{\mathbb{E}}
\newcommand{\md}{\mathbb{D}}
\newcommand{\mn}{\mathbb{N}}

\newcommand{\mr}{\mathbb{R}}
\newcommand{\mm}{\mathcal{M}}
\newcommand{\di}{\widetilde{W}_{n+1}-W_n+R_n}

\newcommand{\od}{\overset{d}{=}}

\begin{document}

\title{On the rate of convergence of a regular martingale related to the
branching random walk}\date{}
\author{A.M. Iksanov\footnote {National T. Shevchenko University of Kiev}\footnote{The original, Ukrainian variant of the paper was published in Ukrainian Mathematical Journal (2006), 58(3), 326-342.} \ \ \ \ \ \ \ \ \ \ \ \
\ \ \ \ \ \ \ \ \ \ \ \ \ \ \ \ \ \ \ \ \ \ \ \  UDK 519.21
}\maketitle
\begin{quote}
Let $\mm_n, n=0,1,\ldots$ be the supercritical branching random
walk, in which the number of direct descendants of one individual
may be infinite with positive probability. Assume that the
standard martingale $W_n$ related to $\mm_n$ is regular, and $W$
is a limit random variable. Let $a(x)$ be a nonnegative function
which regularly varies at infinity, with exponent greater than
$-1$. The paper presents sufficient conditions of the almost sure
convergence of the series $\sum_{n=1}^{\infty}a(n)(W-W_n)$. Also
we establish a criterion of finiteness of $\me W\log^+ W
a(\log^+W)$ and $\me \log^+|\zi| a(\log^+|\zi|)$, where
$\zi:=Q_1+\sum_{n=2}^\infty M_1\cdots M_n Q_{n+1}$, and $(M_n,
Q_n)$ are independent identically distributed random vectors, not
necessarily related to $\mm_n$.
\medskip

\end{quote}

\section{Introduction and main results}

Let $\mm$ be a point process on $\mr$, i.e. random locally finite
counting measure. Assume that $\mm\{+\infty\}=0$ and set $L:=\mm
(\mr)$. In this paper the variable $L$ may be deterministic or
random, finite or infinite with positive probability.

By \emph{the branching random walk} (BRW) we mean the sequence of
point processes $\mm_n, n=0,1,\ldots$, where for any Borel set
$B\subset \mr$, $\mm_0(B)=1_{\{0\in B\}}$,
\begin{equation}\label{BRW}
\mm_{n+1}(B):=\sum_r \mm_{n,r}(B-A_{n,r}), n=0,1,\ldots.
\end{equation}
Here $\{A_{n,r}\}$ are the points of $\mm_n$, and $\{\mm_{n,r}\}$
are independent copies of the $\mm$. More detailed definition of
the process can be found in \cite{IksRos, Iks2004}.

Notice that our definition of the BRW differs from the two
previously known ones. The modern definition of the BRW introduced
in \cite{Big1} assumes that $L<\infty$ a.s. Before the appearance
of \cite{Big1} by the BRW was meant the sequence (\ref{BRW}) where
an underlying point process $\mm$ had independent and identically
distributed points. Now the latter processes are sometimes called
the homogeneous BRW.

In the paper we only consider the \emph{supercritical} BRW.
Therefore, if $\mmp\{L<\infty\}=1$, it is additionally assumed
that $\me L>1$. Recall that supercriticality ensures the survival
of a population with positive probability.

Let $\mathcal{U}:=\bigcup_{n=0}^\infty \mn^n$ be the set of all
finite sequences $u=i_1\ldots i_n, i_k\in \mn$ that contains the
empty sequence $\mn^0:=\{\oslash\}$. A tree $\mathcal{T}$ with
root $\oslash$ is a subset of $\mathcal{U}$ that contains
$\oslash$ and such that $i_1\ldots i_n \in \mathcal{T}$ implies
$i_1\ldots i_k\in \mathcal{T}, k=\overline{1,n-1}$; each element
$i_1\ldots i_n\in \mathcal{T}$ is assigned $L_{i_1\ldots i_n}\in
[0,\infty]$, and $i_1\ldots i_nj\in \mathcal{T} \Leftrightarrow
j\in \{1, \ldots, L_{i_1\ldots i_n}\}$. A tree $\mathcal{T}$ is
called labelled if each $u\in \mathcal{T}$ is assigned a label
$A_u$.

For each realization of the BRW there is a labelled tree with root
$\oslash$. The elements $u$ of this tree are called individuals,
$\oslash$--the initial ancestor; the label $A_u$ defines the
position of the individual $u$ on the real line, $A_\oslash=0$. If
$u=i_1\ldots i_n$, then $n$ is called the generation of the
individual $u$ (notation: $|u|=n$ ($|\oslash|=0$)).

Assume that for some $\gamma>0$
\begin{equation}\label{gam}
m(\gamma):=\int_{-\infty}^\infty e^{\gamma x}\mm_1(dx)\in
(0,\infty).
\end{equation}
For $n=1,2,\ldots$ denote by $\mathcal{F}_n:=\sigma(\mm_1,\ldots,
\mm_n)$ the $\sigma-$field generated by the point processes
$\mm_1, \ldots, \mm_n$ and set
$$W_n:=m^{-n}(\gamma)\int_\mr e^{\gamma x}\mm_n(dx)=m^{-n}(\gamma)\sum_{|u|=n}e^{\gamma A_u}.$$

Under extra moment restrictions in \cite{Big1} and \cite{Liu} (for
the case $L<\infty$ a.s.) and in \cite{Lyons} conditions were
given for the regularity (uniform integrability) of the
non-negative martingale $\{(W_n, \mathcal{F}_n): n=1,2,\ldots\}$.
For the case when $L$ can be infinite with positive probability,
and without extra moment assumptions a criterion of regularity of
the martingale was pointed out in Proposition 1.1 \cite{IksRos}
(see \cite{Iks2004} for a proof).

Recall that the regularity of arbitrary martingale $(U_n,
\mathcal{G}_n)$ ensures the existence of (equivalence class of)
$\mathcal{G}_\infty$--measurable random variable $U$  such that
(a) $\me U=\me U_n$; (b) as $n\to\infty$ $U_n$ a.s. converges to
$U$.

Denote by $W$ the limit random variable for the regular martingale
$W_n$. Then $\me W=1$, and
\begin{equation*}
W= m(\gamma)^{-n}\sum_{|u|=n}e^{\gamma A_u}W^{(u)},
\end{equation*}
where, given $\mathcal{F}_n$, $\{W^{(u)}: |u|=n\}$ are
conditionally independent copies of the $W$.

Set $Y_u:=e^{\gamma A_u}/m^{|u|}(\gamma)$. Let $(Z,S)$ be a random
vector whose distribution is defined by the equality
\begin{equation}
\label{eq}\me \sum_{|u|=1}Y_uk(Y_u, \sum_{|v|=1}Y_v)=\me k(Z,S),
\end{equation}
which is assumed to hold for any nonnegative bounded Borel
function $k(x,y)$.  For problems that the present paper is aimed
at, a joint distribution of $(Z,S)$ does not matter, but knowledge
of marginal distributions is essential. If $k$ does not depend on
$x$, (\ref{eq}) implies the equality
$$\mmp\{S\in dy\}=y\mmp\{W_1\in dy\}.$$ Taking in (\ref{eq})
$k(x,y)=r(x)$ leads to $$\me r(Z)=\me \sum_{|u|=1}Y_ur(Y_u),$$ or,
more generally,
\begin{equation}
\label{eq1} \me r(Z_1\cdots Z_n)=\me \sum_{|u|=n}Y_ur(Y_u),
\end{equation}
where $Z_1, Z_2, \ldots$ are independent copies of the $Z$. Notice
that (\ref{eq1}) holds for any Borel function $r$ with such a
convention: if the right-hand side is infinite or does not exist
the same is true for the left-hand side.

Let $a:\mr^+\to \mr^+$ be a function that regularly varies at
$\infty$ with exponent $\alpha>-1$. If $\alpha=0$, we additionally
assume that $a$ does not decrease near $\infty$. If (\ref{simple})
given below holds, $W_n$ converges to $W$ in mean (see Proposition
1.1 \cite{IksRos}). The paper provides sufficient conditions for
the a.s. convergence of the series
\begin{equation}\label{series}
\sum_{n=0}^\infty a(n)(W-W_n)
\end{equation}
provided (\ref{simple}) holds. This result is a statement about
the rate of the a.s. convergence of the regular martingale $W_n$
to its limit $W$.
\begin{theorem}\label{rateconv}
Let
\begin{equation}\label{simple} \me \log Z \in (-\infty, 0) \
\ \text{and} \ \me W_1 \log^+ W_1<\infty,
\end{equation}
and the distribution of $\log Z$ is non-arithmetic. The conditions
\begin{equation} \label{cond}
\me (\log^+ Z)^3a(\log^+ Z)<\infty \text{ \ \ and \ \ } \me
W_1(\log^+ W_1)^2 a(\log^+ W_1)<\infty
\end{equation}
are sufficient for the a.s. convergence of (\ref{series}).
\end{theorem}
The author thinks that the first inequality in (\ref{cond}) can be
weakened to \newline $\me (\log^+ Z)^2 a(\log^+ Z)<\infty$. If the
conjecture is correct then according to Theorem \ref{mom44} the
following equivalence should be true:
\begin{equation}\label{racon}
\left|\sum_{n=0}^\infty a(n)(W-W_n)\right|<\infty \ \ \text{м.н.}
\Leftrightarrow \me W \log^+ W a(\log^+ W)<\infty.
\end{equation}
The corollary given below proves the conjecture for two particular
cases.
\begin{corr}\label{ra}
Assume that (\ref{simple}) holds. If either $\mathcal{M}(-\infty,
-\gamma^{-1}\log m(\gamma))=0$ a.s., and the distribution of $\log
Z$ is non-arithmetic, or $W_n=\mathcal{M}_n(\mr)/(\me
\mathcal{M}(\mr))^n$, then (\ref{series}) a.s. converges iff $\me
W \log^+ W a(\log^+ W)<\infty$.
\end{corr}
\begin{theorem}\label{mom44}
If (\ref{simple}) holds then
$$\me W\log^+W a(\log^+W)<\infty \text{ \ \ iff \ \ }$$
$$\me W_1 (\log^+W_1)^2 a(\log^+W_1)<\infty.$$
\end{theorem}
\begin{rem}
Theorem 1.3(b)\cite{IksRos} provides a criterion of finiteness of
$\me Wf(W)$ for concave functions $f$ that grow more rapidly than
any power of logarithm. Under the conditions of Theorem
\ref{mom44} this result cannot be used.
\end{rem}

\begin{myproof}\emph{of Theorem \ref{rateconv}}.
We will use the idea of the proof of Theorem 4.1 \cite{Asm}.

On the set of extinction of the population (its probability
measure may equal zero) the series (\ref{series}) contains finite
number of non-zero terms, and hence it trivially converges.
Therefore, not indicating this explicitly, we will investigate the
convergence of the series on the set of survival and assume that
$W>0$.

Without loss of generality we can and do assume that
$m(\gamma)=1$. Indeed, the positions $\{A_u, |u|=n\}$ of
individuals in the $n$-th generation can be replaced by
$\{B_u:=A_u-|u|\log m(\gamma), |u|=n\}$. However we keep all the
previously introduced notation.

Put $b(x):=xa(x)$ and note that $b(x)$ regularly varies at
$\infty$ with exponent $\beta:=\alpha+1>0$. For $n=0,1,\ldots$ we
define the sequences
$$\widetilde{W}_{n+1}:=\sum_{|u|=n}e^{\gamma
A_u}W_1^{(u)}1_{\{b(n)e^{\gamma A_u}W_1^{(u)}\leq 1\}};$$
$$R_n:=\me(W_n-\widetilde{W}_{n+1}|\mathcal{F}_n)=\me
\left(\sum_{|u|=n}e^{\gamma A_u}W_1^{(u)}1_{\{b(n)e^{\gamma
A_u}W_1^{(u)}>1\}}|\mathcal{F}_n\right),$$ where, given
$\mathcal{F}_n$, $\{W_1^{(u)}: |u|=n\}$ are conditionally
independent copies of the random variable $W_1$.
\begin{lemma}\label{ser1}
Assume that (\ref{simple}) holds, and the distribution of $\log Z$
is non-arithmetic. Then the conditions
\begin{equation}
\label{cond100}\me (\log^+ Z)^3 a(\log^+ Z)<\infty \text{ \ \ and
\ \ } \me W_1 \log^+ W_1 a(\log^+W_1)<\infty
\end{equation}
are sufficient for convergence of the series
$$\sum_{n=0}^{\infty}\mmp\{W_{n+1}\neq \widetilde{W}_{n+1}\};$$
$$\sum_{n=0}^{\infty}\mathbb{D}(b(n)(\widetilde{W}_{n+1}-W_n+R_n)).$$
\end{lemma}
Thus, if (\ref{cond100}) holds, the sequence
$$\sum_{n=0}^{m}b(n)(\widetilde{W}_{n+1}-W_n+R_n), m=0,1,\ldots$$
is an $L_2-$ bounded and hence regular martingale. Therefore, the
series
\newline $\sum_{n=0}^{\infty}b(n)(\widetilde{W}_{n+1}-W_n+R_n)$ is
a.s. convergent. By Lemma \ref{ser1} and Borel-Cantelli lemma the
series $\sum_{n=0}^{\infty}b(n)(W_{n+1}-W_n+R_n)$ is a.s.
convergent too. From (4.8) \cite{Asm} it follows that the series
$\sun a(n)(W-W_n+\sum_{k=n}^{\infty}R_k)$ converges a.s.

Consequently, the a.s. convergence of $\sun a(n)(W-W_n)$ is
equivalent to that of $\sun a(n)\sum_{k=n}^{\infty}R_k$, which in
its turn is equivalent to the a.s. convergence of the series $\sun
b(n)R_n$. The latter follows from that fact that $R_n\geq 0$ a.s.,
the equality
$$\sum_{n=1}^m a_n\sum_{k=n}^{\infty}R_k=\left(\sum_{k=1}^m
a_k\right)\sum_{n=m+1}^{\infty}R_n
+\sum_{n=1}^mR_n\left(\sum_{k=1}^n a_k\right),$$ which holds for
any $m \in \mn$, and Lemma 4.2\cite{Asm}.

The next lemma completes the proof of Theorem \ref{rateconv}.
\end{myproof}
\begin{lemma}\label{ser2}
Assume that the conditions (\ref{simple}) and (\ref{cond100})
hold, and the distribution of $\log M$ is non-arithmetic. the
series $\sun b(n)R_n$ is a.s. convergent iff
\begin{equation} \label{part}
\me W_1(\log^+ W_1)^2 a(\log^+W_1)<\infty.
\end{equation}
\end{lemma}

At this point it is appropriate to prove Corollary \ref{ra}.

\begin{myproof}\emph{of Corollary \ref{ra}}.
Let $\mathcal{M}(-\infty, -\gamma^{-1}\log m(\gamma))=0$ a.s., or
equivalently $Z\in [0,1]$ a.s. In this case, the inequalities
containing $Z$ in Theorem \ref{rateconv} and Theorem \ref{mom44},
and Lemma \ref{ser1} and Lemma \ref{ser2} hold automatically.
Assume that $\me W \log^+ W a(\log^+ W)<\infty$. By Theorem
\ref{mom44} this is equivalent to
$$\me W_1 (\log^+ W_1)^2 a(\log^+W_1)<\infty.$$ By Theorem \ref{rateconv} the series (\ref{series}) is a.s. convergent.

Let now the series (\ref{series}) is a.s. convergent. If $\me
W_1\log^+ W_1 a(\log ^+W_1)<\infty$, Lemma \ref{ser2} implies that
$\me W_1(\log^+ W_1)^2 a(\log ^+W_1)<\infty$. Therefore, in view
of Theorem \ref{mom44} $\me W \log^+ W a(\log ^+W) <\infty$.
Assume that\newline $\me W_1\log^+ W_1 a(\log ^+W_1)=\infty$.
Since by the assumption of the corollary $\me
W_1\log^+W_1<\infty$, then $\alpha\geq 0$. If $\alpha>0$, then
there exists a $\delta\in [0,\alpha)$ such that $\me W_1(\log^+
W_1)^{\delta+1}<\infty$ and $\me W_1(\log^+
W_1)^{\delta+2}=\infty$. By Lemma \ref{ser2} the series $\sun
n^\delta (W-W_n)$ diverges. According to Abel's criterion, for
$\epsilon\in (0, \alpha-\delta)$ the series $\sun
n^{\alpha-\epsilon}(W-W_n)$ cannot converge. Hence the series
(\ref{series}) diverges. If $\alpha=0$, then $\me
W_1\log^+W_1<\infty$ and $\me W_1(\log^+W_1)^2=\infty$. By Lemma
\ref{ser2} the series $\sun (W-W_n)$ diverges. By the assumption
at the beginning of the section $a(x)$ does not decrease for large
$x$. This implies that the series (\ref{series}) diverges. The
proof of the corollary for Galton-Watson process follows a similar
route. It suffices to remark that in this case we should take
$\gamma=0$ in (\ref{gam}) and that $Z=(\me \mathcal{M}(\mr))^{-1}$
a.s.
\end{myproof}

\begin{myproof}\emph{of Lemma \ref{ser1}}.
Denote by $F(x)$ the distribution function of the random variable
$W_1$. Let $S_n$ be a random walk, starting at the origin, with a
step distributed like $(-\log Z)$. By the assumption of the lemma
$\mu:=\me S_1\in (0,\infty)$. By Lemma \ref{pert}(b) for $x>0$
\begin{equation}
\label{V} V(x):=\sun b(n)\mmp\{S_n \leq \log b(n)+\log x\}<\infty.
\end{equation}
For $x>0$ define
$$K(x):=\int_0^x ydV(y)=xV(x)-\int_0^x V(y)dy;$$
$$M(x):=\int_x^{\infty}y^{-1}dV(y)=-x^{-1}V(x)+\int_x^{\infty}y^{-2}V(y)dy.$$ Since the function
$l(x):=\mu^{-\alpha-2} b(\log x)$ slowly varies at $\infty$, and
by Lemma \ref{pert} $V(x)$ defined in (\ref{V}) satisfies
(\ref{Haan}), this $V$ belongs to de Haan's class $\Pi_l$.

By Theorem 3.7.1 \cite{BGT}
\begin{equation}
\label{imp} \lix \dfrac{K(x)}{xb(\log x)}=\mu^{\alpha+2}; \lix
\dfrac{xM(x)}{b(\log x)}=\mu^{\alpha+2}.
\end{equation}
Further we have
$$\sum_{n=0}^\infty \mmp\{W_{n+1}\neq \widetilde{W}_{n+1}\}=
\sum_{n=0}^{\infty}\mmp\{b(n)\underset{|u|=n}{\sup}e^{\gamma
A_u}W_1^{(u)}>1\}\leq$$ $$\leq
\sum_{n=0}^{\infty}\me\sum_{|u|=n}\mmp\{b(n)e^{\gamma
A_u}W_1^{(u)}>1|\mathcal{F}_n\}=$$
$$=\sum_{n=0}^{\infty}\me\sum_{|u|=n}e^{\gamma
A_u}\left(\int_{b^{-1}(n)e^{-\gamma A_u}}^{\infty}dF(x)e^{-\gamma
A_u}\right)\overset{(\ref{eq1})}{=}\sum_{n=0}^{\infty}\me
e^{S_n}\int_{b^{-1}(n)e^{S_n}}^{\infty}dF(x)=$$
$$=\int_0^{\infty}\me\left(\sum_{n=0}^{\infty}e^{S_n}1_{\{e^{S_n}\leq
b(n)x\}}\right)dF(x)=\int_0^{\infty}K(x)dF(x).$$ The latter
integral converges in view of (\ref{imp}) and $\me W_1
b(\log^+W_1)<\infty$.

Recall that the conditional variance is defined by the equality
$\md(X|\mathcal{G})=\me (X^2|\mathcal{G})-(\me(X|\mathcal{G}))^2$.
Since $\me (\widetilde{W}_{n+1}|\mathcal{F}_n)=W_n-R_n$ and $\me
(\di|\mathcal{F}_n)=0$, then $$\md (b(n)(\di))=b^2(n)\me\left(\md
(\widetilde{W}_{n+1}|\mathcal{F}_n)\right).$$ Further $$\md
(\widetilde{W}_{n+1}|\mathcal{F}_n)=\sum_{|u|=n}\md(e^{\gamma
A_u}W_1^{(u)}1_{\{b(n)e^{\gamma A_u}W_1^{(u)}\leq
1\}}|\mathcal{F}_n)\leq$$ $$\leq \me \left(\sum_{|u|=n}e^{2\gamma
A_u}\me (W_1^21_{\{b(n)e^{\gamma A_u}W_1\leq
1\}})|\mathcal{F}_n\right)=$$
\begin{equation}
\label{disp} =\me \left(\sum_{|u|=n}e^{2\gamma
A_u}\int_0^{b^{-1}(n)e^{-\gamma
A_u}}x^2dF(x)|\mathcal{F}_n\right).
\end{equation}
Thus
$$\su \md (b(n)(\di))=\su b^2(n)\me\left(\md
(\widetilde{W}_{n+1}|\mathcal{F}_n)\right)\overset{(\ref{eq1}),
(\ref{disp})}{\leq}$$
$$\leq \su b^2(n)\me e^{-S_n}\int_0^{b^{-1}(n)e^{S_n}}x^2dF(x)=\int_0^{\infty}x^2\me\left(\su b^2(n)e^{-S_n}1_{\{e^{S_n}>
b(n)x\}}\right)dF(x)=$$ $$=\int_0^{\infty}x^2M(x)dF(x).$$ The
latter integral converges in view of (\ref{imp}) and $\me
W_1b(\log^+W_1)<\infty$.
\end{myproof}

For each fixed $x\in\mr$ consider the random variables
$$Q(x):=\sun b(n)\sum_{|u|=n}e^{\gamma A_u}1_{\{e^{\gamma
A_u}>e^{-x}\}};$$
$$\widehat{Q}(x):=\sun
b(n)\sum_{|u|=n}e^{\gamma A_u}1_{\{e^{\gamma
A_u}>e^{-x}b^{-1}(n)\}}.$$ If $$\me \log Z\in (-\infty, 0) \
\text{and} \ \me (\log^+Z)^2 b(\log^+Z)<\infty,$$ these are a.s.
finite. This follows from Theorem 1 \cite{Als} that guarantees
that for all $x\in \mr$
$$\me Q(x)=\sun b(n)\mmp\{S_n\leq x\}<\infty,$$
where $S_n$ is the same random walk as in the proof of Lemma
\ref{ser1}. That $\me \widehat{Q}(x)<\infty$ follows from similar
considerations and inequality (\ref{cond112}).
\begin{lemma}\label{ser3}
If (\ref{simple}) holds, and $\me (\log^+Z)^2 b(\log^+Z)<\infty$,
than a.s. on the set of survival
\begin{equation}
\label{ren} \lix \dfrac{Q(x)}{xb(x)}= \lix
\dfrac{\widehat{Q}(x)}{xb(x)} =\dfrac{W}{(\beta+1)(-\me \log
Z)^{\beta+1}}>0,
\end{equation}
where $\beta>0$ is the exponent of regular variation of $b$.
\end{lemma}
\begin{myproof} is similar to the proof of Theorem B
\cite{Big1}. Pick any $0<a<\mu=-\me\log Z$. For each $x>0$ there
exists an integer $N=N(x)>0$ such that $(N-1)^2a\leq x<N^2 a$. For
$x>0$ define the random variables
$$Q_1(x):=\dfrac{1}{(N-1)^2ab((N-1)^2a)}\left(\sum_{n=1}^{N^2}b(n)W_n\right),$$
$$Q_2(N,x):=\dfrac{1}{(N-1)^2ab((N-1)^2a)}\left(\sum_{n=N^2}^{\infty}b(n)
\sum_{|u|=n}e^{\gamma A_u}1_{\{e^{\gamma A_u}>e^{-an}\}}\right).$$
Recall that for almost all $\omega$ from the set of survival
$W(\omega)>0$. Since, as $m\to\infty$, $\sum_{n=1}^{m}b(n)W_n\sim
W\sum_{n=1}^{m}b(n)$ a.s.; $\sum_{n=1}^{m}b(n)\sim
(\beta+1)^{-1}mb(m)$, then
\begin{equation}
\label{1}\lix Q_1(x)= \dfrac{W}{(\beta+1)a^{\beta+1}} \text{ \ \
a.s.   \ \ }
\end{equation}
By Lemma \ref{pert}(а) the series $\sun b(n)\mmp\{S_n-an\leq 0\}$
converges. Hence $$\me
\sum_{N=2}^{\infty}Q_2(N,x)=\sum_{N=2}^{\infty}\dfrac{1}{(N-1)^2ab((N-1)^2a)}\left(\sum_{n=N^2}^{\infty}b(n)
\mmp\{e^{-S_n}>e^{-an}\}\right)<\infty,$$ which implies that
\begin{equation}
\label{2} \lix Q_2(N,x)=0  \text{ \ \ a.s.  \ \ }
\end{equation}
By Theorem 1.5.3 \cite{BGT} without loss of generality we can and
do assume that $b(x)$ does not decrease for $x>0$. Therefore for
$x>0$ $\dfrac{Q(x)}{xb(x)}\leq Q_1(x)+Q_2(x)$. Now taking into
account (\ref{1}) and (\ref{2}), and sending $a\to \mu$ gives
\begin{equation}
\label{inter} \lisx \dfrac{Q(x)}{xb(x)}\leq
\dfrac{W}{(\beta+1)\mu^{\beta+1}} \text{ \ \ a.s. \ \ }
\end{equation}

Pick now any $a>\mu$. For each $x \geq a$ there exists an integer
$N=N(x)>0$ such that $Na\leq x<(N+1)a$. For $x\geq a$ consider the
random variables
$$Q_3(x):=\dfrac{1}{(N+1)ab((N+1)a)}\left(\sum_{n=1}^Nb(n)W_n\right),$$
$$Q_4(x):=\dfrac{1}{(N+1)ab((N+1)a)}\left(\sum_{n=0}^N b(n)
\sum_{|u|=n}e^{\gamma A_u}1_{\{e^{\gamma A_u}\leq
e^{-an}\}}\right).$$ In a similar way as it was done for $Q_1(x)$
we can prove that
\begin{equation}
\label{3}\lix Q_3(x)= \dfrac{W}{(\beta+1)a^{\beta+1}} \text{ \ \
a.s. \ \ }
\end{equation}
The proof of the fact that
\begin{equation}
\label{4}\lix Q_4(x)=0 \text{ \ \ a.s. \ \ }
\end{equation}
almost coincides with the proof of a similar statement given on p.
35 \cite{Big1}. By Theorem 4.2 \cite{Spitzer} $r:=\sun n^{-1}
\mmp\{S_n>an\}<\infty$. Therefore
$$\me \sun n^{-1}\sum_{|u|=n}e^{\gamma A_u}1_{\{e^{\gamma A_u}\leq
e^{-an}\}}=r<\infty.$$ By Kroneker's lemma
$$\lin (n+1)^{-1}\sum_{k=1}^n \sum_{|u|=n}e^{\gamma
A_u}1_{\{e^{\gamma A_u}\leq e^{-an}\}}=0 \text{ \ \ a.s. \ \ }$$
From this and monotonicity of $b$ (\ref{4}) follows.

For large $x$ $\dfrac{Q(x)}{xb(x)}\geq Q_3(x)-Q_4(x)$. Hence from
(\ref{3}) and (\ref{4}) letting $a$ go to $\mu$ we get
$$\liix \dfrac{Q(x)}{xb(x)}\geq \dfrac{W}{(\beta+1)\mu^{\beta+1}}
\text{ \ \ a.s. \ \ }$$ Together with (\ref{inter}) the latter
inequality proves the limit relation for $Q(x)$.

Now fix $\delta \in (0,\mu)$ and choose $r=r(\delta)>0$ so that
$\log b(n)\leq \delta n+r, n=1,2,\ldots$ We have the following
\begin{equation}\label{cond112} Q(x)\leq
\widehat{Q}(x)\leq \sun b(n)\sum_{|u|=n}e^{\gamma
A_u}1_{\{e^{\gamma A_u}>e^{-x-\delta n-r}\}}.
\end{equation}
Using the same analysis as above allows us to check that the
right-hand side of this inequality satisfies the same limit
relation (\ref{ren}) as $Q(x)$.
\end{myproof}

\begin{myproof}\emph{of Lemma \ref{ser2}}. The definition of
$R_n$ implies the following representation
$$R_n=\sum_{|u|=n}e^{\gamma A_u}\int_{b^{-1}(n)e^{-\gamma
A_u}}^\infty xdF(x),$$ where, as before, $F(x)$ is the
distribution function of $W_1$. We have the (formal) equality
\begin{gather*}
\label{5}\sun b_n R_n=\int_0^\infty x dF(x)\sun
b(n)\sum_{|u|=n}e^{\gamma A_u}1_{\{\gamma A_u>-\log x-\log
b(n)\}}=
\\ =\int_0^\infty \widehat{Q}(\log x)x dF(x).
\end{gather*}
By the assumptions of the lemma $\me \log Z \in (-\infty, 0)$ and
$\me (\log^+Z)^2 b(\log^+Z)<\infty$. Therefore by Lemma
\ref{ser3}, as $x\to\infty$, $\widehat{Q}(\log x)\sim const \log x
b(\log x)$ a.s. Hence, the series with non-negative terms $\sun
b(n)R_n$ converges iff (\ref{part}) holds.
\end{myproof}

\section{Moments of random series and the proof of Theorem \ref{mom44}}
Assume that $(M_1, Q_1),(M_2, Q_2), \ldots$ are independent copies
of a random vector $(M, Q)$, not necessarily related to the BRW,
and defined on a fixed probability space. Set
$$\Pi_0:=1 \ \ \text{and}  \ \ \Pi_n:=M_1 M_2 \cdots M_n,
n=1,2,\ldots;$$
\begin{equation}\label{per} \zi:=\sum_{k=1}^\infty
\Pi_{k-1}Q_k
\end{equation}
Throughout this section we assume that
\begin{equation*}
\mmp\{M=0\}=0, \mmp\{Q=0\}<1,
\end{equation*}
and, provided the series in (\ref{per}) is absolutely convergent,
that the distribution of $\zi$ is non-degenerate.

Theorem \ref{mom4} given below may be of some interest on its own,
supplements the result of Theorem 1.6 \cite{IksRos} and is a key
ingredient in proving Theorem \ref{mom44}. Recall that $b(x)$ is a
regularly varying function with exponent $\beta>0$ and put
$c(x):=xb(x)$.
\begin{theorem}\label{mom4} If
\begin{equation}\label{gol}
\me \log |M| \in (-\infty, 0) \ \text{and} \ \me \log^+|Q|<\infty,
\end{equation}
then
\begin{equation}
\label{int}\me b(\log^+|\zi|)<\infty \Leftrightarrow \me c(\log^+
|M|)<\infty, \ \ \me c(\log^+ |Q|)<\infty.
\end{equation}
\end{theorem}
\begin{proof} The functions $b$ and $c$ can be represented as follows: $b(x)=x^\beta
L(x)$, $c(x)=x^{\beta+1}L(x)$, where $L(x)$ slowly varies at
$\infty$. For $y>1$ set
$\Lambda_{\beta}(y):=\dfrac{\log^{\beta-1}yL(\log y)}{\beta y}$.
This function regularly varies at $\infty$ with exponent $(-1)$.
The function $\underset{t\geq x}{\sup}\Lambda_{\beta}(t)$ does not
increase, and by Theorem  1.5.3 \cite{BGT} $\underset{t\geq
x}{\sup}\Lambda_\beta(t)\sim \Lambda_\beta(x)$ (from here on the
record $F\sim G$ means that $\lix (F(x)/G(x))=1$). Changing
variable and then appealing to Karamata's theorem yield
\begin{equation}
\label{rel} b(\log x)\sim \int_1^x \Lambda_\beta (y)dy \sim
\int_1^x \underset{t\geq y}{\sup}\Lambda_\beta (t)dy=:
\tilde{f}(x-1).
\end{equation}
Analogously
\begin{equation*}
c(\log x)\sim \int_1^x\Lambda_{\beta+1}(y)dy \sim \int_1^x
\underset{t\geq y}{\sup}\Lambda_{\beta+1}(t)dy=:\phi(x-1).
\end{equation*}
The functions $\tilde{f}$ and $\phi$ are non-decreasing and
concave on $\mr^+$. Also these equal $0$ at $x=0$ and tend to
$\infty$ as $x\to\infty$. In particular, $\phi$ is subadditive.
Also, from (\ref{rel}) and Karamata's theorem it follows that
\begin{equation*}
(\beta+1)^{-1} c(\log x)\sim \int_1^x
(\tilde{f}(y)/y)dy=:\tilde{g}(x-1).
\end{equation*}
The function $c(x)$ regularly varies at $\infty$ with exponent
$\beta+1>1$. By Theorem 1.5.3 \cite{BGT} it is equivalent at
$\infty$ to a non-decreasing function. Therefore, according to
Lemma 1(a) \cite{Als} there exists a non-decreasing function
$\psi(x)\sim c(\log x)$ such that $\psi(x)=0$ for $x\leq 1$ and
\begin{equation}
\label{ineq}\psi(xy)\leq a(\psi(x)+\psi(y))
\end{equation}
for all $x,y\in [1,\infty)$ and some positive constant $a$.

Thus we conclude that it suffices to prove the equivalence
(\ref{int}) with $b(\log x)$ replaced by $\tilde{f}(x)$, and
$c(\log x)$--by $\tilde{g}(x)$, $\phi(x)$ or $\psi(x)$.

To prove implication $\Leftarrow$ of the theorem we should use the
fact that according to Theorem 2.1 \cite{GolMal} the condition
(\ref{gol}) ensures that $|\zi|<\infty$ a.s.

First we assume that $|M|\in [0,1]$ a.s. In that case the
condition \newline $\me c(\log^+|M|)<\infty$ holds automatically.
Let $\me c(\log^+|Q|)<\infty$ or equivalently $\me
\tilde{g}(|Q|)<\infty$. By Theorem 1.6(a)\cite{IksRos} $\me
\tilde{f}(|\zi|)<\infty$. This is equivalent to $\me
b(\log^+|\zi|)<\infty$. The proof of the other way implication
goes the same path and, in particular, appeals to the same Theorem
1.6(a)\cite{IksRos}.

Now we are ready to consider the general case. First assume that
in (\ref{int}) the inequalities for $|M|$ and $|Q|$ hold.
Equivalently,
\begin{equation*}
\me \tilde{g}(|M|)<\infty, \me \tilde{g}(|Q|)<\infty.
\end{equation*}

Consider the random variables
$$N_0:=0, N_{i+1}:=\inf\{n>N_{i}: |\Pi_n|<|\Pi_{N_i}|\},
i=0,1,\ldots$$ Since under the assumptions of the theorem
$\Pi_n\to 0$ a.s. as $n\to\infty$, $\me N_i<\infty, i=1,2,\ldots$.
For $k=1,2,\ldots$ put
\begin{equation*}
M_k^\prime:=|M_{N_{k-1}+1}|\cdots|M_{N_{k}}|, \Pi_0^\prime:=1,
\Pi_k^\prime:=M_1^\prime \cdots M_k^\prime;
\end{equation*}
\begin{equation*}
Q_k^\prime:=|Q_{N_{k-1}+1}|+|M_{N_{k-1}+1}||Q_{N_{k-1}+2}|+\cdots+|M_{N_{k-1}+1}|
\cdots|M_{N_{k}-1}||Q_{N_{k}}|.
\end{equation*}
The random vectors $\{(M_k^\prime, Q_k^\prime), k=1,2,\ldots\}$
are independent copies of \newline $(|\Pi_{N_1}|,
\sum_{k=1}^{N_1}|\Pi_{k-1}||Q_k|)$ and, furthermore,
\begin{equation*}
\sum_{k=1}^\infty
|\Pi_{k-1}||Q_k|=\sum_{k=1}^{\infty}\Pi_{k-1}^\prime Q_k^\prime.
\end{equation*}
If we could prove that
\begin{equation}
\label{in10}\me
\tilde{g}(\sum_{k=1}^{N_1}|\Pi_{k-1}||Q_k|)<\infty,
\end{equation}
than this implied that $\me b(|\zi|)<\infty$, and, therefore, one
way of the theorem would be established. Indeed, since
$|\Pi_{N_1}|\in (0,1)$ a.s. and $\mmp\{|\Pi_{N_1}|=1\}=0$, and
(\ref{in10}) guarantees that $\me \log^+
(\sum_{k=1}^{N_1}|\Pi_{k-1}||Q_k|)<\infty$, than the first part of
the proof applied on the vector $(|\Pi_{N_1}|,
\sum_{k=1}^{N_1}|\Pi_{k-1}||Q_k|)$ instead of $(|M|,|Q|)$ gives
the wanted.

Let us check (\ref{in10}) with $\tilde{g}$ replaced by $\psi$.
Since
$$\sum_{k=1}^{N_1}|\Pi_{k-1}||Q_k| \leq N_1
\underset{1\leq k\leq N_1}{\sup}|\Pi_{k-1}||Q_k| \leq N_1
\underset{0\leq k\leq N_1-1}{\sup}|\Pi_k| \sum_{i=1}^{N_1}
|Q_i|,$$ then taking into account (\ref{ineq}) allows us to
conclude that to prove (\ref{in10}) it suffices to establish three
inequalities: 1) $\me \psi(N_1)<\infty$; \newline 2) $\me
\psi(\underset{0\leq k\leq N_1-1}{\sup}|\Pi_k|)<\infty$; 3) $\me
\psi(\sum_{i=1}^{N_1} |Q_i|)<\infty$. Since $\me N_1<\infty$, and
$\psi$ grows more slowly than the linear function, the first
inequality holds true. Further, $\psi(e^x)$ regularly varies with
exponent $\beta+1>1$. Therefore, according to (\ref{defKKK}), the
second inequality is implied by $\me \psi(|M|\vee 1)<\infty$. The
latter is equivalent to $\me c(\log^+|M|)<\infty$. To check the
third inequality we replace $\psi$ with $\phi$. A benefit of the
replacement is that $\phi$ be subadditive. As the random variables
$1_{\{N_1 \geq n\}}$ and $|Q_n|$ are independent, we have
$$\me \phi(\sum_{i=1}^{N_1} |Q_i|)\leq \me \sum_{i=1}^{N_1}
\phi(|Q_i|)=\me N_1 \me \phi(|Q|)<\infty.$$

Now assume that the left-hand side of (\ref{int}) holds. This is
equivalent to
\begin{equation}\label{cond11}
\me \tilde{f}(|\zi|)<\infty.
\end{equation}
By Proposition 3.1 \cite{IksRos}, either $\infty>\me
\tilde{f}(\underset{n\geq 0}{\sup}|\Pi_n|)$ or $\infty>\me
\tilde{f}(\underset{n\geq 0}{\sup}|\Pi_{2n}|)$. Equivalently,
either $\infty>\me f(\underset{n\geq 0}{\sup} S_n)$, or
$\infty>\me f(\underset{n\geq 0}{\sup} \grave{S}_n)$, where
$S_n:=\log |\Pi_n|$, $\grave{S}_n:=\log |\Pi_{2n}|, n=0,1,\ldots$
are random walks with steps distributed like $\log |M|$ and $\log
|M_1M_2|$ respectively. In view of (\ref{defKKK}), either $\me
g(\log^+M)<\infty$ or $\me g(\log^+(M_1M_2))<\infty$. Clearly,
both of these imply the next to last inequality.

On the other hand, by Proposition 3.1 \cite{IksRos} (\ref{cond11})
implies that either
$$\me \tilde{f}(\underset{k\geq
1}{\sup}\Pi^\ast_{k-1}|Q_k^s|)\leq \me \tilde{f}(\underset{k\geq
1}{\sup}|\Pi_{k-1}||Q_k^s|)<\infty, \ \ \text{or}$$ $$\me
\tilde{f}(\underset{k\geq
1}{\sup}\grave{\Pi}^\ast_{k-1}|\grave{Q}^s_k|)\leq \me
\tilde{f}(\underset{k\geq
1}{\sup}|\grave{\Pi}_{k-1}||\grave{Q}^s_k|)<\infty,
$$ hold, where $$\grave{\Pi}_0:=1, \ \
\grave{\Pi}_n:=\grave{M}_1 \grave{M}_2\cdots \grave{M}_n,
n=1,2,\ldots,$$ the vectors $$(\grave{M}_k,
\grave{Q}_k):=(M_{2k-1}M_{2k}, M_{2k-1}Q_{2k}+Q_{2k-1}), \ \
k=1,2,\ldots$$ are independent and identically distributed; $(M_n,
Q_n)\od (M_n, Q_n^\prime)$, given $M_n$, $Q_n$ and $Q_n^\prime$
are conditionally independent, $Q_n^s:=Q_n-Q_n^\prime$,
$\grave{Q}_n^s$ and $\grave{Q}_n^\prime$ have the same meaning,
but are defined in terms of $\grave{M}_n$ and $\grave{Q}_n$;
$\Pi_0^\ast:=1$, $\Pi_k^\ast:=M_1^\ast\cdots M_k^\ast$,
$M_k^\ast:=|M_k|\wedge 1, k=1,2,\ldots$, and $\grave{\Pi}_k^\ast$
are defined similarly. Since $M_k^\ast, \grave{M}_k^\ast \leq 1$
a.s., and strictly smaller than one with positive probability,
Corollary 3.1 \cite{IksRos} implies that $\me
\tilde{g}(|Q|)<\infty$. Hence, $\me g(\log^+|Q|)<\infty$.
\end{proof}

\begin{myproof}\emph{of Theorem \ref{mom44}}. Theorem \ref{mom44}
can be obtained from Theorem \ref{mom4} by using the same approach
that was exploited in \cite{IksRos} to deduce Theorem 1.3 from
Theorem 1.6. Theorem \ref{mom4} applies to the random series
generated by the vector $(Z,S)$. The latter was defined in
(\ref{eq}).
\end{myproof}

\section{Appendix}

Lemma \ref{pert} is a key ingredient in proving Lemma \ref{ser1}.
Part b) of the lemma deals with \emph{perturbed} random walks and
generalizes a result of \cite{Als} for random walks.
\begin{lemma}\label{pert}
Assume that a function $\varphi:\mr^+\to\mr^+$ regularly varies at
$\infty$ with exponent $\beta>0$. Let $T_n, n=0,1,\ldots$ be a
random walk starting at zero with $\mu:=\me T_1\in (0,\infty)$. If
$\me (T_1^-)^2\varphi(T_1^-)<\infty$, then\newline a) for any
$\epsilon>0$
$$I:=\sun \varphi(n)\mmp\{T_n>(\mu+\epsilon)n\}<\infty;{   } \sun \varphi(n)\mmp\{T_n\leq
(\mu-\epsilon)n\}<\infty;$$ b) for all $x\in \mr$
$$V(x):=\sun \varphi(n)\mmp\{T_n\leq \log \varphi(n)+\log x\}<\infty,$$ if, moreover, the distribution of
$T_1$ is non-arithmetic, then for all $h>0$
\begin{equation}
\label{Haan} \lix \dfrac{V(hx)-V(x)}{\mu^{-\beta-1}\varphi(\log
x)}=\log h.
\end{equation}
\end{lemma}
\begin{proof} According to Theorem 1.5.3 \cite{BGT} we can and do assume that $\varphi$ is non-decreasing on $\mr^+$.\newline a) The
sequence $\widetilde{T}_n:=-T_n+(\mu+\epsilon)n, n=0,1,\ldots$ is
a random walk with $\me \widetilde{T}_1=\epsilon \in (0,\infty)$.
Therefore by Theorem 1(a) \cite{Als} $I=\sun
\varphi(n)\mmp\{\widetilde{T}_n\leq 0\}<\infty$. The convergence
of the second series can be established similarly. \newline (b) In
what follows we will use the idea of the proof of Theorem 2
\cite{Lai}. Fix $\delta \in (0, \mu)$ and pick $r=r(\delta)>0$
such that $\log \varphi(n)\leq \delta n+r, n=1,2,\ldots$. The
sequence $\widehat{T}_n:=T_n-\delta n$ is a random walk with $\me
\widehat{T}_1=\mu-\delta \in (0,\infty)$. Since $V(x)\leq \sun
\varphi(n)\mmp\{\widehat{T}_n\leq \log x+r\}$, and the latter
series converges by Theorem 1(a)\cite{Als}, then the function
$V(x)$ is finite for all $x>0$. The relation (\ref{Haan}) is
tantamount to
\begin{equation}
\label{Haan1} \lix
\dfrac{U(x+h)-U(x)}{\mu^{-\beta-1}\varphi(x)}=h, \text{ \ \ for
all \ \ } h\in \mr,
\end{equation}
where $U(x):=V(e^x).$ Actually it suffices to prove (\ref{Haan1})
for small positive $h$ from an interval $(h_0, h_1)$ (see, for
example, Lemma 3.2.1 \cite{BGT}). Fix such an $h$.
\newline For any $\epsilon\in (0, \mu/2)$ and large enough $x$
the inequality $\log \varphi(n)\leq \epsilon n$ holds for $n\geq
N_2=N_2(x):=\left[\dfrac{x+h}{\mu-2\epsilon}+1\right]$. Set
$N_1=N_1(x):=\left[\dfrac{x+h}{\mu+\epsilon}\right]$. By using
part a) of the lemma, given $\rho>0$ we can choose $m=m(\rho)>0$
such that
\begin{equation}
\label{ad}\sum_{n=m+1}^\infty
\varphi(n)\mmp\{T_n>(\mu+\epsilon)n\}\leq \rho.
\end{equation}
Write $$U(x+h)-U(x)=\sun \varphi(n)\mmp\{T_n\leq \log
\varphi(n)+x\}=$$$$=\sum_{n=1}^{m}+\sum_{n=m+1}^{N_1}+\sum_{n=N_1+1}^{N_2-1}+\sum_{n=N_2}^{\infty}=:
I_1(x)+I_2(x)+I_3(x)+I_4(x).$$ It is obvious that $\lix
I_1(x)=0$.\newline If for large $x$ and $n\geq N_2(x)$ $T_n\geq
(\mu-\epsilon)n$, then $T_n-\log \varphi(n)\geq (\mu-2\epsilon)n$.
Therefore, when $x\to \infty$, $$I_4(x)\leq
\sum_{n=N_2(x)}^{\infty}\varphi (n)\mmp\{T_n\leq
(\mu-\epsilon)n\}\to 0$$ by part a) of the lemma.\newline If for
large $x$ and $n\in \{m+1,\ldots, N_1(x)\}$ $T_n\leq
(\mu+\epsilon)n$, then $$T_n-\log\varphi(n)\leq
(\mu+\epsilon)N_1-h\leq x.$$ Hence $$I_2(x)\leq
\sum_{n=m+1}^{N_1}\varphi(n)\mmp\{T_n-\log \varphi(n)>x\}\leq$$
$$\leq \sum_{n=m+1}^{N_1}\varphi(n)\mmp\{T_n>(\mu+\epsilon)n\}\overset{(\ref{ad})}{\leq}
\rho.$$ By Potter's inequality (Theorem 1.5.6 \cite{BGT}), for any
positive $q$ and $\theta$ there exists an $x_0>0$ such that
$$\log\varphi (\dfrac{x+h}{\mu-2\epsilon})-\log\varphi(\dfrac{x+h}{\mu+\epsilon})\leq (1+q)+(\beta+\theta)
(\log(\mu+\epsilon)-\log(\mu-2\epsilon)):=B(q,\theta).$$ Hence for
$x\geq x_0$
$$I_3(x)\leq \sum_{n=N_1+1}^{N_2-1}\varphi(n)\mmp\{\log\varphi(N_1+1)+x<T_n\leq
\log\varphi(N_2-1)+x+h\}\leq$$
$$\sun\varphi(n)\mmp\{\log \varphi(\dfrac{x+h}{\mu-2\epsilon})+x<T_n\leq \log\varphi(\dfrac{x+h}{\mu-2\epsilon})+x+h+B(q.\theta)\}.$$
An appeal to Theorem 2 \cite{Als} yields
$$\underset{x\to\infty}{\lim\sup}\dfrac{I_3(x)}{\varphi(x)}\leq
\dfrac{h+B(q,\theta)}{\mu^{\beta+1}}.$$ Letting $q$ and $\epsilon$
go to $0$ results in
$$\underset{x\to\infty}{\lim\sup}\dfrac{I_3(x)}{\varphi(x)}\leq
\dfrac{h}{\mu^{\beta+1}}.$$ Thus, we have proved that
$$\underset{x\to\infty}{\lim\sup}\dfrac{U(x+h)-U(x)}{\varphi(x)}\leq
\dfrac{h}{\mu^{\beta+1}}.$$ We now intend to check that
\begin{equation}
\label{low}
\underset{x\to\infty}{\lim\inf}\dfrac{U(x+h)-U(x)}{\varphi(x)}\geq
\dfrac{h}{\mu^{\beta+1}}.
\end{equation}
Put $R_n:=T_n-\log \varphi(n), n=1,2,\ldots$. For each $\epsilon
\in (0,\mu)$ we define
$N_3=N_3(x):=\left[\dfrac{x+h}{\mu-\epsilon}+1\right]$ and make
use of the random variable $N_1$ defined above. For any positive
$q$ and $\theta$ such that $\tau=\tau(q,\theta, \epsilon):=\log
(1+q)(\dfrac{\mu+\epsilon}{\mu-\epsilon})^{\beta+\theta}<h_0$, and
large $x$ Potter's inequality $\log \varphi(N_3(x)-1)-\log \varphi
(N_1(x))\leq \tau$ holds. Moreover, we have the following
$$U(x+h)-U(x)\geq \sum_{n=N_1+1}^{N_3-1}\varphi(n)\mmp\{x<R_n\leq
x+h\}\geq$$
$$\geq\sum_{n=N_1+1}^{N_3-1}\varphi(n)\mmp\{\log\varphi(n)-\log\varphi(N_1)+x<R_{N_1}+T_n-T_{N_1}\leq
x+h\}\geq$$ $$\geq
\sum_{n=1}^{N_3-N_1-1}\varphi(n+N_1)\mmp\{\tau+x-R_{N_1}<T_n\leq
x-R_{N_1}+h\}\geq$$ $$\geq
\varphi(\dfrac{x}{\mu+\epsilon})\sum_{n=1}^{N_3-N_1-1}\mmp\{\tau+x-R_{N_1}<T_n\leq
x-R_{N_1}+h\}=$$
$$=\varphi(\dfrac{x}{\mu+\epsilon})\me g(x-R_{N_1(x)}),$$ where $g(t):=\sum_{n=1}^{N_3-N_1-1}\mmp\{\tau+t<T_n\leq
t+h\}$. We will show that a.s.
\begin{equation}
\label{g}\lix g(x-R_{N_1(x)})=\mu^{-1}(h-\tau).
\end{equation}
Blackwell's theorem \cite{Bl} implies that the function $g(t)$ is
bounded. Therefore from (\ref{g}) it follows that $$\lix \me
g(x-R_{N_1(x)})=\mu^{-1}(h-\tau).$$ Consequently, taking into
account the regular variation of $\varphi$ allows us to conclude
that
$$\underset{x\to\infty}{\lim\inf}\dfrac{U(x+h)-U(x)}{\varphi(x)}\geq
\dfrac{h-\tau(q,\theta, \epsilon)}{(\mu+\epsilon)^{\beta}\mu}.$$
Sending $q$ and $\epsilon$ to $0$ leads to (\ref{low}).

By the strong law of large numbers, as $x\to\infty$ we have
$R_{N_1(x)}=\mu N_1(x)+o(N_1(x))$ a.s. Hence, as $x\to\infty$
$x-R_{N_1(x)}=\epsilon(\mu+\epsilon)^{-1}x+o(x)$ a.s. To prove
(\ref{g}) it suffices to verify that for \emph{arbitrary}
non-random function $z(x)=\epsilon(\mu+\epsilon)^{-1}x+o(x)$
\begin{equation}
\label{inte} \lix
\sum_{n=1}^{N_2(x)-N_1(x)-1}\mmp\{\tau+z(x)<T_n\leq
z(x)+h\}=\mu^{-1}(h-\tau).
\end{equation}
If positive integer $n\geq N_3-N_1$ and $T_n>(\mu-\epsilon)n$,
then for large $x$ $T_n>2\epsilon x(\mu+\epsilon)^{-1}+h>z(x)+h$.
Therefore
$$\sum_{n=N_3(x)-N_1(x)}^{\infty}\mmp\{T_n\leq z(x)+h\}\leq \sum_{n=N_3(x)-N_1(x)}^{\infty}\mmp\{T_n\leq (\mu-\epsilon)n\}.$$
According to part a) of the lemma the last expression tends to $0$
when $x\to\infty$. By Blackwell's theorem (\ref{inte}) holds and
hence (\ref{g}) holds too.
\end{proof}

Let $\xi_1, \xi_2, \ldots$ be independent copies of a random
variable $\xi$ with $m:=\me \xi \in (-\infty, 0)$. Define
$S_0:=0$, $S_n:=\xi_1+\ldots+\xi_n, n=1,2,\ldots$. Then
$M_\infty:=\suin S_n<\infty$ a.s., and $\me \tau^-_x<\infty$,
where
$$\tau^-_x:=\inf\{n: S_n<-x\} \ , \ x\geq 0 .$$
Let $f$ be a non-negative measurable function such that
$\underset{x\to \infty}{\lim}f(x)=\infty$ and there exists an
$x_0\geq 0$ such that $f$ is increasing and concave for $x\geq
x_0$. Define the new function $g$ as follows:
\begin{equation*}
g(x):=\int_{x_0}^x(f(y)/y)dy \ \ \text{for} \ \ x\geq x_0; \ \
g(x):=0 \ \ \text{for} \ \ x<x_0.
\end{equation*}
Set $u(x):=f(e^x)$, $v(x):=g(e^x)$. Assume that a function $h$
regularly varies at $\infty$ with exponent $\beta>0$.
\begin{lemma}\label{Alsm}
For $x\geq 0$
\begin{equation} \label{defK}
\me u(M_\infty)<\infty \Leftrightarrow \me v(\underset{0\leq n\leq
\tau^-_x-1}{\sup}S_n)<\infty.
\end{equation}
Each of these inequalities ensures that
\begin{equation} \label{defKK}
\me v(\xi^+)<\infty.
\end{equation}
Also the following equivalences hold:
\begin{equation}\label{defKKK} (\underset{0\leq n\leq \tau^-
-1}{\sup}S_n)h(\underset{0\leq n\leq \tau^--1}{\sup}S_n)<\infty
\Leftrightarrow  \me h(M_\infty)<\infty \Leftrightarrow \me \xi^+
h(\xi^+)<\infty.
\end{equation}
\end{lemma}
\begin{proof}
Without loss of generality we can assume that $f$ is increasing
and concave on $\mr^+$, $f(0)=0$, $\lix f(x)=\infty$, and
$g(x)=\int_0^x(f(u)/u)du$. This follows from the fact that we can
consider the function $\widehat{f}(x)=f(x+x_0)-f(x_0)$ in place of
$f$. This function possesses the properties listed above, and the
ratio $\widehat{f}/f$ is bounded away from zero and bounded from
the above. For fixed $x\geq 0$ define the random variables
$N_0:=0$,
$$N_{i+1}:=\inf\{n>N_i: S_n<S_{N_i}-x\}, i=0,1,\ldots.$$ Notice that $\tau^-_x=N_1$ and $N_i<\infty$ a.s. Put
$$V_k:=\sup\{S_{N_k}, S_{N_k+1},\ldots, S_{N_{k+1}-1}\}, k=0,1,\ldots;$$
$$Z_{k+1}:=\sup \{0, \xi_{N_k+1}, \ldots, \xi_{N_k+1}+\ldots+\xi_{N_{k+1}-1}\} ,k=0,1,\ldots.$$
Then $V_k=S_{N_k}+Z_{k+1}$ and $M_\infty=\underset{k\geq
0}{\sup}V_{k}$. Notice that $Z_1, Z_2,\ldots$ are independent
copies of $Z:=\underset{0\leq i \leq \tau^-_x-1}{\sup}S_i$, and by
Wald's identity $\me |Z|\leq \me \sum_{k=1}^{\tau^-_x}|\xi_k|=\me
\tau^-_x \me |\xi|<\infty$.\newline $\Leftarrow$ in (\ref{defK}).
Since $S_{N_k}<-kx$, and for fixed $\epsilon>0$
$$\mmp\{M_\infty>y\}\leq \mmp\{\underset{k\geq 0}{\sup}(-kx+Z_{k+1})>y)\leq $$
$$\leq \sum_{k=0}^\infty \mmp\{Z_{k+1}>y+k(x+\epsilon)\}\leq
\sum_{k=[y/(x+\epsilon)]}^\infty \mmp\{Z>k(x+\epsilon)\}\leq$$
$$\leq \int_{[y/(x+\epsilon)]-1}^\infty \mmp\{Z>(x+\epsilon)y\}dy.$$
The latter integral converges, since $\me |Z|<\infty$. Thus,
$$\infty>\me u(M_\infty)=\int_0^\infty u^\prime(z)\mmp\{M_\infty>z\}dz,$$ if
$$\int_{-\infty}^\infty u^\prime(z)\int_z^\infty \mmp\{Z>y\}dydz<\infty.$$
Since $u(x)=v^\prime (x)$, integrating by parts shows that the
latter inequality is equivalent to
$$\infty>\int_{-\infty}^\infty v^\prime (z)\mmp\{Z>z\}dz=\me v(Z)=\me v(\underset{0\leq n\leq
\tau^-_x-1}{\sup}S_n).$$
\newline $\Rightarrow$ in (\ref{defK}). $\{S_{N_{k}}, k=1,2,\ldots\}$ is a random walk starting at zero and with a step
distributed like $S_{\tau^-_x}$. The random vectors
$(S_{N_k}-S_{N_{k-1}}, Z_k), k=1,2,\ldots$ are independent and
identically distributed, and $\lin S_{N_n}=-\infty$ a.s. Let
$(\widetilde{M}_1, \widetilde{Q}_1),(\widetilde{M}_2,
\widetilde{Q}_2), \ldots $ be independent copies of the vector
$(\widetilde{M}:=e^{S_{\tau^-_x}}, \widetilde{Q}:=e^Z)$. By
construction $\mmp\{\widetilde{M}\leq 1\}=1$ and
$\mmp\{\widetilde{M}=1\}=0$. Therefore, by Corollary 3.1
\cite{IksRos} the inequality $\me g(\widetilde{Q})<\infty$ follows
from $\me f(\suik \widetilde{M}_1\cdots \widetilde{M}_{k-1}
\widetilde{Q}_k)<\infty$. It remains to note that
$\widetilde{Q}=\exp(Z)=\exp(\underset{0\leq i \leq
\tau^-_x-1}{\sup}S_i)$. In a similar way $\suik
\widetilde{M}_1\cdots
\widetilde{M}_{k-1}\widetilde{Q}_k=\exp(\suk(S_{N_{k}}+Z_{k+1}))=\exp(M_\infty)$.

Since $\xi_1^+\leq \underset{0\leq n\leq \tau^-_x-1}{\sup}S_n$,
(\ref{defKK}) is implied by (\ref{defK}).

At the beginning of the proof of Theorem \ref{mom4} it is shown
that there exists a non-decreasing, concave on $\mr^+$ function
$f$ that additionally satisfies $f(0)=0$, $\lix f(x)=\infty$ and
$h(x)\sim f(e^x)$. Therefore the first equivalence and implication
$\Leftarrow$ in the second equivalence in (\ref{defKKK}) follow
from (\ref{defK}). The rest can be deduced from Theorem 3
\cite{Als}.
\end{proof}

\end{document}